\documentclass{article} 
\usepackage{url}

\usepackage{amsmath, amsthm, amssymb}
\usepackage{mathtools, bm}
\usepackage{manfnt}
\usepackage{xfrac}
\usepackage{graphicx}
\usepackage[inline]{enumitem}
\usepackage[colorlinks=true, allcolors=blue]{hyperref}
\usepackage[svgnames]{xcolor}
\newcommand*{\greysquare}{\textcolor{gray}{\blacksquare}}
\usepackage{cleveref}

\DeclareMathOperator{\tr}{tr}
\DeclareMathOperator{\adj}{adj}

\begin{document}
%
\begin{center}
    \Large Parametric-ROM of Structures with Varying Geometry using Direct Parameterization of Invariant Manifolds
\end{center}
%
%
\begin{center}
    Tiago Martins\textsuperscript{1,2},
    Alessandra Vizzaccaro\textsuperscript{3} and 
    Daniel Rixen\textsuperscript{2}
\end{center}

%
%
%

\noindent
\textsuperscript{1} tiago.martins@tum.de\\
\textsuperscript{2}
Chair of Applied Mechanics, 
TUM School of Engineering and Design,\\
Department of Mechanical Engineering,
Technical University of Munich\\
Boltzmannstr. 15 
D - 85748 Garching, Germany\\
\textsuperscript{3}
College of Engineering, Mathematics and Physical
Sciences\\
University of Exeter, Exeter, UK

\begin{abstract}
This work presents a framework for parametric reduction in FEM, where geometry is controlled by a parameter without altering material properties or stress states. The inverse determinant in the weak form is expanded as a power series, with explicit expressions for the zeroth and first-order terms. External forcing and parameter dependence are incorporated into an enlarged autonomous system, reduced via the direct parameterization of invariant manifolds method and homological equations. The parameter is treated as an additional variable with trivial dynamics, isolated for inclusion in the ROM. This approach enables efficient parametric studies and advances reduced-order modeling in structural dynamics.

\end{abstract}

\section{Introduction}

Cabré et al.~\cite{cabre2003parameterization1,cabre2003parameterization2,cabre2005parameterization3} introduced the parametrization method, later refined by Haro et al.~\cite{haro2016parameterization} to enhance reduction to invariant manifolds. While earlier methods focused on invariant manifolds or normal form theory~\cite{guckenheimer2013nonlinear,touze2004hardening}, the direct parametrization revealed both could stem from the solvable invariance equation using graph or normal-form styles. Researchers later adapted this to vibratory systems with dissipation~\cite{haller2016nonlinear}, automating reduction for two-dimensional manifolds with damping and nonlinearities. 
The development of MORFE2.0~\cite{vizzaccaro2024direct} and SSMtool 2.0~\cite{jain2022compute} automated high-order SSM approximations for FEM. Opreni et al.~\cite{opreni2023high}, Grolet et al.~\cite{grolet2024high} and Bettini et al.~\cite{bettini2024model} further extended these methods to accommodate generic forcing, parametric excitation, non-polynomial nonlinearities, and piecewise smooth dynamics.

In earlier work, Marconi et al.~\cite{marconi2021higher} developed a parametric nonlinear reduced-order model for structures with geometric imperfections. They split the displacement into two components: one due to parametric defects and another due to structural motion, expanding each independently. The method describes defects through user-defined displacement fields and expresses the internal forces as polynomials in both the defect and displacement fields. The Neumann expansion approximates the strains and simplifies the elastic forces into polynomials.

Most recently, Morsy and and Tiso~\cite{morsy2025predicting} used polynomial chaos expansion to analyze how uncertainties such as ISO tolerances, surface imperfections, and random bolt preloads affect the dynamics of bolted joints. They evaluated variability in nonlinear modal responses with minimal simulations, showing how parameter changes impact damping and stiffness, making it a clear parametric approach to uncertainty quantification.

Building on these methods, this paper contributes to the parametric reduction framework by expanding as a power series the inverse of the determinant of the geometry morphing.

\section{Equation of motion of solid mechanics}
Consider a body with volume \mbox{$V$} and surface \mbox{$S$}, defined by material properties: density \mbox{$\rho$} and elasticity tensor \mbox{$\mathcal{A}$}. Let \mbox{$\mathbf{u}$} be the displacement field, \mbox{$\mathbf{b}$} the body force, and \mbox{$\mathbf{t}$} the surface traction. In weak-form the equation of motion becomes
\begin{gather}
\begin{split}
    &\int_{V} \rho \,  \ddot{\mathbf{u}} \cdot \mathbf{w} + (\nabla \mathbf{u})_\mathrm{sym} : \mathcal{A} : (\nabla \mathbf{w})_\mathrm{sym} \, \mathrm{d}V \\
    + &\int_{V} (\nabla \mathbf{u})_\mathrm{sym} : \mathcal{A} : (\nabla^{\mathrm{T}} \mathbf{u} \nabla \mathbf{w})_\mathrm{sym} 
    + \frac{1}{2} \nabla^{\mathrm{T}} \mathbf{u} \nabla \mathbf{u} : \mathcal{A} : (\nabla \mathbf{w})_\mathrm{sym} \, \mathrm{d}V\\
    + &\int_{V} \frac{1}{2} \nabla^{\mathrm{T}} \mathbf{u} \nabla \mathbf{u} : \mathcal{A} : (\nabla^{\mathrm{T}} \mathbf{u} \nabla \mathbf{w})_\mathrm{sym} \, \mathrm{d}V\\
    = &\int_{V} \mathbf{b} \cdot \mathbf{w} \, \mathrm{d}V + \int_{S} \mathbf{t} \cdot \mathbf{w} \, \mathrm{d}S \, ,
\end{split}
\label{eq:EOM}
\end{gather}
where $(\greysquare)_\mathrm{sym}$ is the symmetric part of a matrix. This should be satisfied for any suficiently continuous test field \mbox{$\mathbf{w}$} vanishing on Dirichlet boundaries of \mbox{$\mathbf{u}$}. The expression above separates the different orders by rows:
\begin{enumerate*}
    \item writes the terms that are linear on the displacement and its derivatives,
    \item the quadratic terms,
    \item the cubic, and
    \item the external excitation, independent on the displacement.
\end{enumerate*}
Moreover, including Rayleigh damping is straightforward after this development.

\section{Structures with parametric geometry}
The objective is to perform parametric reduction when the geometry of the body is controlled via a parameter. In FEM, one controls nodal positions, thereby modifying the shape of each element and the body. 

\subsection{How to parameterize geometry}

\noindent
\textbf{Approach 1: Deformation with relaxation}\\
Deform the body and reset internal stress (relaxation). This process reshapes the structure while conserving total mass, leading to density redistribution. The parameter controls the deformation and the relaxation forces that eliminate pre-stress.

\vspace{0.5\baselineskip}\noindent
\textbf{Approach 2: Deformation with pre-stress}\\
Deform the body while preserving internal stress (no relaxation), and conserve total mass by redistributing density. The parameter characterizes the external pre-stress force, which alters the system's equilibrium configuration. Alternatively, the parameter defines the pre-stressed deformation state.

\vspace{0.5\baselineskip}\noindent
\textbf{Approach 3: Modify the undeformed geometry}\\
Reshape the geometry of the body and preserve material properties. In FEM, the meshes of the varying configurations must be compatible via a continuously differentiable map. Namely, the parameter controls node positions without modifying density, the stress state, or introducing additional forces.

\vspace{\baselineskip}
This work addresses Approach 3, which is more challenging as the parameter alters the integration domain in the weak-form of the equation of motion (EOM). The first two approaches, that entail additional parametric-forces in the EOM, are more straightforward to handle.

\subsection{Modifying the undeformed geometry: mathematical formulation}
\sloppy The quantities in the original configuration are marked with the subscript $\greysquare_0$, and $\nabla_{\! 0}$ is the derivative w.r.t. $\bm{x}_0$. Given the geometry parameter $\mu$, the altered configuration $V$ is described by the map \mbox{$\bm{x} = \bm{x}_0 + \bm{\mho}_1(\bm{x}_0) \mu$}. Here, the analysis is restricted to orientation-preserving maps satisfying \mbox{$\det(\nabla_{\! 0} \bm{x})>0$}, such that the volume differential becomes $\mathrm{d}V = \det(\nabla_{\! 0} \bm{x}) \, \mathrm{d}V_0$. The deformation gradient w.r.t $\bm{x}$ can be expressed as $\nabla \mathbf{u} = \nabla_{\! 0} \mathbf{u} \, (\nabla_{\! 0} \bm{x})^{-1}$. Additionally, the determinant and the adjugate can expressed as polynomials on the parameter $\mu$. Expressly,
\begin{align*}
    \det(\nabla_{\! 0} \bm{x}) = 1 + \tr(\nabla_{\! 0}\bm{\mho}_1) \, \mu + \tr(\adj(\nabla_{\! 0}\bm{\mho}_1)) \, \mu^2 + \det(\nabla_{\! 0}\bm{\mho}_1) \, \mu^3 \, , \\
    \adj(\nabla_{\! 0} \bm{x}) = (\nabla_{\! 0} \bm{x})^{-1} \, \det(\nabla_{\! 0} \bm{x}) = \mathbf{I} + [\tr(\nabla_{\! 0}\bm{\mho}_1) \mathbf{I} - \nabla_{\! 0}\bm{\mho}_1] \, \mu + \adj(\nabla_{\! 0}\bm{\mho}_1) \, \mu^2 \, .
\end{align*}

\begin{proof}
    The determinant is formally related to the characteristic polynomial as follows: \mbox{$\det(\mathbf{I} + \bm{B} \mu) = -\mu^3 \det(\lambda \, \mathbf{I} - \bm{B})$} with \mbox{$\lambda \, \mu = -1$}.
    Regarding the adjugate matrix \mbox{$\bm{C}= \adj(\mathbf{I} + \bm{B} \mu)$}, one has, per definition
    \begin{align}\label{eq:adjoint_definiton}
        (\mathbf{I} + \bm{B} \mu) \, \bm{C} = \det(\mathbf{I} + \bm{B} \mu) \, \mathbf{I} \, .
    \end{align}
     The degrees w.r.t. $\mu$ are \mbox{$\deg_\mu (\mathbf{I} + \bm{B} \mu)\le 1$} and \mbox{$\deg_\mu (\det(\mathbf{I} + \bm{B} \mu))\le 3$}, therefore \mbox{$\deg_\mu (\bm{C})\le 2$}. Matching the orders of $\mu$ on equation \eqref{eq:adjoint_definiton} yields
    \begin{align*}
        &\text{order $\mu^0$}: & \bm{C}_0 = \mathbf{I} &\implies \bm{C}_0 = \mathbf{I}\\
        &\text{order $\mu^1$}: & \bm{C}_1 + \bm{B} \, \bm{C}_0 = \tr(\bm{B}) \, \mathbf{I} &\implies \bm{C}_1 = \tr(\bm{B}) \mathbf{I} - \bm{B}\\
        &\textcolor{lightgray}{\text{order $\mu^2$}:} & \textcolor{lightgray}{\bm{C}_2 + \bm{B} \, \bm{C}_1 = \tr(\adj(\bm{B})) \, \mathbf{I}} & \textcolor{lightgray}{\phantom{\implies} \; \; \text{redundant equation}}\\
        &\text{order $\mu^3$}: & \bm{B} \, \bm{C}_2 = \det(\bm{B}) \, \mathbf{I} &\implies \bm{C}_2 = \adj(\bm{B}) \, . \qed\\
    \end{align*}
\end{proof}


When the integration domain of the weak-form is changed to $V_0$, the term \mbox{$(\nabla_{\! 0} \bm{x})^{-1} = \adj(\nabla_{\! 0} \bm{x}) \det(\nabla_{\! 0} \bm{x})^{-1}$} emerges from the deformation gradient. The adjugate is a polynomial on $\mu$, but the inverse of the determinant remains problematic. This work suggests two methods.

\vspace{\baselineskip}\noindent
\textbf{Method 1: Additional scalar field}\\
One introduces a scalar field \mbox{$c: V_0 \to \mathbb{R}$} constrained by \mbox{$c \cdot \det(\nabla_{\! 0} \bm{x}) = 1$}. With test fields $w(\bm{x}_0)$, this can be transformed into weak-form as
\begin{align}
    \int_{V_0} w \, \Big[ c \cdot \det(\nabla_{\! 0} \bm{x}) - 1 \Big] \, \mathrm{d} V_0 = 0 \, .
\end{align}
In practice, $c$ must be discretized, which introduces additional variables with accompanying algebraic equations. The amount of degrees of freedom of the resulting DAE is one-third larger, when using the same mesh as in the FEM. Moreover, the accuracy depends on the choice of shape functions and quadrature method. Previous work~\cite{frangi2023reduced} used a field $\mathbf{c}$ constrained by \mbox{$\mathbf{c} \nabla_{\! 0} \bm{x} = \mathbf{I}$}.

\vspace{\baselineskip}\noindent
\textbf{Method 2: Power series}\\
Considering \mbox{$\det(\nabla_{\! 0} \bm{x}) = 1 + h_1 \mu + h_2 \mu^2 + h_3 \mu^3$}, express $\det(\nabla_{\! 0} \bm{x})^{-1}$ as a power series in $\mu$, given that \mbox{$0<\det(\nabla_{\! 0} \bm{x})<2$}. Namely, 
\begin{align}
    &\det(\nabla_{\! 0} \bm{x})^{-1} = \frac{1}{1 - [1-\det(\nabla_{\! 0} \bm{x})]} = \sum_{n=0}^{\infty} (1-\det(\nabla_{\! 0} \bm{x}))^n = \sum_{\kappa=0}^\infty a_\kappa \mu^\kappa
\end{align}
where the polynomial coefficients are as follows
\begin{align*}
    a_\kappa = \sum_{\bm{\sigma} \in S_\kappa} \frac{(\sigma_1+\sigma_2+\sigma_3)!}{\sigma_1! \; \sigma_2! \; \sigma_3!} (-h_1)^{\sigma_1} (-h_2)^{\sigma_2} (-h_3)^{\sigma_3}  \, ,
\end{align*}
with \mbox{$S_k := \left\{\bm{\sigma} \in \mathbb{N}_{0}^3 : \sigma_1 + 2 \sigma_2 + 3 \sigma_3 = \kappa \right\}$}.

\vspace{0.5\baselineskip}
In previous work~\cite{marconi2021higher}, $(\nabla_{\! 0} \bm{x})^{-1}$ was directly expanded as a power series in the parameter $\mu$, with the coefficients being \mbox{$3\times 3$} matrices. Within their radius of convergence, these power series expansions achieve greater accuracy as additional terms are incorporated. However, analyzing excessively high orders of the parameter $\mu$ proves impractical since the invariant manifold reduction method already imposes a finite truncation order in its asymptotic expansion.

\vspace{\baselineskip}
Increasing the amount of degrees of freedom undermines the efficiency of reduced order modeling. Therefore, this work utilizes the power series method because it preserves the number of DoFs in the FOM. Consequently, the equation of motion splits into multiple equations, one at every order of the parameter.

\subsection{Equation of motion for varying orders of the parameter}
The equations of motion (EOM) decompose into terms of different powers of $\mu$. At order $\mu^0$, the EOM remains in its standard form. At order $\mu^1$, the EOM includes additional terms that capture the linear dependence of the dynamics on $\mu$. At higher orders of $\mu$, the expressions grow progressively more intricate.

\vspace{\baselineskip}\noindent
\textbf{Order $\mu^0$}\\
The equations of motion simplify to the standard governing equations of solid mechanics, expressed in their weak-form as:
\begin{gather*}
\begin{split}
    &\int_{V_0} \rho \,  \ddot{\mathbf{u}} \cdot \mathbf{w} \, \mathrm{d}V_0 + \int_{V_0} (\nabla_{\! 0} \mathbf{u})_\mathrm{sym} : \mathcal{A} : (\nabla_{\! 0} \mathbf{w})_\mathrm{sym} \, \mathrm{d}V_0 \\
    + &\int_{V_0} (\nabla_{\! 0} \mathbf{u})_\mathrm{sym} : \mathcal{A} : (\nabla_{\! 0}^{\mathrm{T}} \mathbf{u} \nabla_{\! 0} \mathbf{w})_\mathrm{sym} 
    + \frac{1}{2} \nabla_{\! 0}^{\mathrm{T}} \mathbf{u} \nabla_{\! 0} \mathbf{u} : \mathcal{A} : (\nabla_{\! 0} \mathbf{w})_\mathrm{sym} \, \mathrm{d}V_0\\
    + &\int_{V_0} \frac{1}{2} \nabla_{\! 0}^{\mathrm{T}} \mathbf{u} \nabla_{\! 0} \mathbf{u} : \mathcal{A} : (\nabla_{\! 0}^{\mathrm{T}} \mathbf{u} \nabla_{\! 0} \mathbf{w})_\mathrm{sym} \, \mathrm{d}V_0\\
    = &\int_{V_0} \mathbf{b} \cdot \mathbf{w} \, \mathrm{d}V_0 + \int_{S_0} \mathbf{t} \cdot \mathbf{w} \, \mathrm{d}S_0 \, ,
\end{split}
\end{gather*}

\vspace{\baselineskip}\noindent
\textbf{Order $\mu^1$}\\
The linear terms, i.e. the inertia and the linear stiffness forces, are given by
\begin{gather*}
\begin{split}
    \int_{V_0} &\Big[\rho \,  \ddot{\mathbf{u}} \cdot \mathbf{w} + (\nabla_{\! 0} \mathbf{u})_\mathrm{sym} : \mathcal{A} : (\nabla_{\! 0} \mathbf{w})_\mathrm{sym} \Big] \tr(\nabla_{\! 0} \bm{\mho}_1) \\
    - &\Big[ (\nabla_{\! 0} \mathbf{u})_\mathrm{sym} : \mathcal{A} : (\nabla_{\! 0} \mathbf{w} \, \nabla_{\! 0} \bm{\mho}_1)_\mathrm{sym} + (\nabla_{\! 0} \mathbf{u} \, \nabla_{\! 0} \bm{\mho}_1)_\mathrm{sym} : \mathcal{A} : (\nabla_{\! 0} \mathbf{w})_\mathrm{sym} \Big] \, \mathrm{d}V_0 \, .\\ 
\end{split}
\end{gather*}
The quadratic stiffness term expands to
\begin{gather*}
\begin{split}
    \int_{V_0} &\Big[ (\nabla_{\! 0} \mathbf{u})_\mathrm{sym} : \mathcal{A} : (\nabla_{\! 0}^{\mathrm{T}} \mathbf{u} \nabla_{\! 0} \mathbf{w})_\mathrm{sym} 
    + \frac{1}{2} \nabla_{\! 0}^{\mathrm{T}} \mathbf{u} \nabla_{\! 0} \mathbf{u} : \mathcal{A} : (\nabla_{\! 0} \mathbf{w})_\mathrm{sym} \Big] \tr(\nabla_{\! 0} \bm{\mho}_1)\\
    - &\Big[ (\nabla_{\! 0} \mathbf{u} \nabla_{\! 0} \bm{\mho}_1)_\mathrm{sym} : \mathcal{A} : (\nabla_{\! 0}^{\mathrm{T}} \mathbf{u} \nabla_{\! 0} \mathbf{w})_\mathrm{sym} 
    + (\nabla_{\! 0}^\mathrm{T} \bm{\mho}_1 \nabla_{\! 0}^{\mathrm{T}} \mathbf{u} \nabla_{\! 0} \mathbf{u})_\mathrm{sym} : \mathcal{A} : (\nabla_{\! 0} \mathbf{w})_\mathrm{sym} \\
    &\phantom{\int_{V_0}\,} + (\nabla_{\! 0} \mathbf{u})_\mathrm{sym} : \mathcal{A} : (\nabla_{\! 0}^\mathrm{T} \bm{\mho}_1 \nabla_{\! 0}^{\mathrm{T}} \mathbf{u} \nabla_{\! 0} \mathbf{w})_\mathrm{sym} + \frac{1}{2} \nabla_{\! 0}^{\mathrm{T}} \mathbf{u} \nabla_{\! 0} \mathbf{u} : \mathcal{A} : (\nabla_{\! 0} \mathbf{w} \nabla_{\! 0}^\mathrm{T} \bm{\mho}_1)_\mathrm{sym}\\
    &\phantom{\int_{V_0}\,} + (\nabla_{\! 0} \mathbf{u})_\mathrm{sym} : \mathcal{A} : (\nabla_{\! 0}^{\mathrm{T}} \mathbf{u} \nabla_{\! 0} \mathbf{w} \nabla_{\! 0} \bm{\mho}_1)_\mathrm{sym} \Big] \, \mathrm{d}V_0 \, .\\
\end{split}
\end{gather*}
The cubic stiffness term becomes
\begin{gather*}
\begin{split}
    \int_{V_0} \frac{1}{2} &\Big[ \nabla_{\! 0}^{\mathrm{T}} \mathbf{u} \nabla_{\! 0} \mathbf{u} : \mathcal{A} : (\nabla_{\! 0}^{\mathrm{T}} \mathbf{u} \nabla_{\! 0} \mathbf{w})_\mathrm{sym} \Big] \tr(\nabla_{\! 0} \bm{\mho}_1) \\
    - &\Big[ (\nabla_{\! 0}^{\mathrm{T}} \mathbf{u} \nabla_{\! 0} \mathbf{u}\nabla_{\! 0} \bm{\mho}_1)_\mathrm{sym} : \mathcal{A} : (\nabla_{\! 0}^{\mathrm{T}} \mathbf{u} \nabla_{\! 0} \mathbf{w})_\mathrm{sym} \\
    &\phantom{\int_{V_0}\,} +\frac{1}{2} \nabla_{\! 0}^{\mathrm{T}} \mathbf{u} \nabla_{\! 0} \mathbf{u} : \mathcal{A} : (\nabla_{\! 0} \bm{\mho}_1^\mathrm{T} \nabla_{\! 0}^{\mathrm{T}} \mathbf{u} \nabla_{\! 0} \mathbf{w})_\mathrm{sym} \\
    &\phantom{\int_{V_0}\,} +\frac{1}{2} \nabla_{\! 0}^{\mathrm{T}} \mathbf{u} \nabla_{\! 0} \mathbf{u} : \mathcal{A} : (\nabla_{\! 0}^{\mathrm{T}} \mathbf{u} \nabla_{\! 0} \mathbf{w} \nabla_{\! 0} \bm{\mho}_1)_\mathrm{sym}
    \Big] \, \mathrm{d}V_0 \, .
\end{split}
\end{gather*}

\section{Direct parameterization of invariant manifolds}
Without loss of generality, the ODEs considered are first-order and autonomous: they take the form \mbox{$\dot{\bm{y}} = \mathbf{A}(\bm{y})$}. When external excitations are introduced, the system becomes non-autonomous (time-dependent); however, constructing an equivalent autonomous representation is achievable by appending auxiliary equations that govern the forcing terms. The state of the external forcing is modeled as a solution of an appropriately defined ODE, \mbox{$\dot{\bm{y}}_\mathrm{ext} = \mathbf{A}_\mathrm{ext}(\bm{y}_\mathrm{ext})$}. The internal dynamics then evolve based on the internal and external states:~\mbox{$\dot{\bm{y}}_\mathrm{int} = \mathbf{A}_\mathrm{int}(\bm{y}_\mathrm{int}, \bm{y}_\mathrm{ext})$}. To incorporate a parameter $\mu$ into the dynamics, the same approach applies, treating the parameter as an additional variable with trivial dynamics, represented by $\dot{\mu} = 0$. To retain the forcing and the parameter in the ROM, one must isolate the corresponding variables and governing equations for inclusion in the ROM, before applying the procedure.

We analyze the full-order dynamical system \mbox{$\dot{\bm{y}} = \mathbf{A}(\bm{y})$} of dimension $n$ and seek to derive a ROM of dimension $m\ll n$ with governing equation \mbox{$\dot{\bm{z}} = \mathbf{f}(\bm{z})$}. The ROM must preserve the essential dynamics of the high-dimensional system near a stable equilibrium point, assumed without loss of generality to be at the origin, satisfying \mbox{$\mathbf{A}(\mathbf{0}) = \mathbf{0}$}. By operating in a reduced-dimensional space, the ROM facilitates computationally efficient simulations and analysis. The ROM trajectories map to the full system’s phase space via the manifold mapping~\mbox{$\bm{y} = \mathbf{W}(\bm{z})$}, ensuring accurate representation of the system's behavior.

This approach assumes that the FOM dynamics $\mathbf{A}(\bm{y})$, the ROM dynamics~$\mathbf{f}(\bm{z})$, and the manifold mapping $\mathbf{W}(\bm{z})$ can all be expanded into power series: 
\begin{align*}
    \mathbf{A}(\bm{y}) = \sum_{\kappa \in \mathbb{N}} \mathbf{A}_\kappa \bm{y}^{\otimes \kappa} \qquad \mathbf{W}(\bm{z}) = \sum_{\kappa \in \mathbb{N}} \mathbf{W}_{\!\kappa} \bm{z}^{\otimes \kappa} \qquad \mathbf{f}(\bm{z}) = \sum_{\kappa \in \mathbb{N}} \mathbf{f}_\kappa \bm{z}^{\otimes \kappa} \, ,
\end{align*}
where $\mathbf{A}_\kappa$ is \mbox{$n\times n^\kappa$}, $\mathbf{W}_\kappa$ is \mbox{$n\times m^\kappa$}, $\mathbf{f}_\kappa$ is \mbox{$m\times m^\kappa$}, and the Kronecker power $\greysquare^{\otimes \kappa}$ represents the repeated Kronecker product: for~$\kappa = 0$, it is the size-$1$ identity matrix; and for $\kappa \ge 1$, it corresponds to $\greysquare \otimes \greysquare \otimes \dotsm \otimes \greysquare$ ($\kappa$ times).
Directly substituting the power series ansatz writes
\begin{align*}
    \nabla_ {\!\bm{z}} \mathbf{W} (\bm{z}) \; \mathbf{f}(\bm{z}) &\approx \mathbf{A}(\mathbf{W}(\bm{z})) \; \iff \\
    \left(\sum_{\kappa \in \mathbb{N}} \mathbf{W}_{\!\kappa} \, \nabla_ {\!\bm{z}} \bm{z}^{\otimes \kappa}\right) \; \left(\sum_{\sigma \in \mathbb{N}} \mathbf{f}_\sigma \bm{z}^{\otimes \sigma} \right) &\approx \sum_{\kappa \in \mathbb{N}} \mathbf{A}_\kappa \left(\sum_{\sigma \in \mathbb{N}} \mathbf{W}_{\!\sigma} \bm{z}^{\otimes \sigma}\right)^{\!\! \otimes \kappa} .
\end{align*}
The term $\nabla_ {\!\bm{z}} \bm{z}^{\otimes \kappa}$ can be expanded into the following order \mbox{$\kappa-1$} expression 
\begin{align*}
    \nabla_ {\!\bm{z}} \bm{z}^{\otimes \kappa} = \sum_{\sigma=1}^{\kappa} \bm{z}^{\otimes \kappa - \sigma} \otimes \mathbf{I}_m \otimes  \bm{z}^{\otimes \sigma - 1} \, .
\end{align*}
To prove this by induction on $\kappa$, use the property of the derivative of a Kronecker product: \mbox{$\nabla_ {\!\bm{z}} (\mathbf{h} \otimes \mathbf{g} ) = \nabla_ {\!\bm{z}} \mathbf{h} \otimes \mathbf{g} + \mathbf{h} \otimes \nabla_ {\!\bm{z}} \mathbf{g}$}.

We aim to equate polynomial terms of the same degree, up to some threshold, and accumulate the error in higher-order terms, similarly to Taylor expansions. 

\vspace{\baselineskip}\noindent
\textbf{Left-hand side}\\
Collecting the terms of order $p$, from the LHS, yields
\begin{align*}
    \text{LHS}\rvert_p = \sum_{\kappa = 1}^{p} \mathbf{W}_{\! \kappa} \Big[\nabla_ {\!\bm{z}} \bm{z}^{\otimes \kappa} \Big] \, \mathbf{f}_{\sigma} \, \bm{z}^{\otimes \sigma} = \sum_{\kappa = 1}^{p} \mathbf{W}_{\! \kappa} \Big[\nabla_ {\!\bm{z}} \bm{z}^{\otimes \kappa} \Big] \, \mathbf{f}_{p - \kappa+1} \, \bm{z}^{\otimes p - \kappa+1} \, ,
\end{align*}
where \mbox{$\sigma$} is such that \mbox{$(\kappa - 1) + \sigma = p$}, given that \mbox{$\kappa - 1$} is the degree of \mbox{$\nabla_ {\!\bm{z}} \bm{z}^{\otimes \kappa}$}. Namely, \mbox{$\sigma = p - \kappa+1$} which bounds \mbox{$1 < \kappa < p$}, since \mbox{$\sigma, \kappa \in \mathbb{N}$}. Subsequently, expanding the derivative yields
\begin{align*}
    &\text{LHS}\rvert_p = \sum_{\kappa = 1}^{p} \sum_{\sigma=1}^{\kappa} \mathbf{W}_{\! \kappa} \Big[ \bm{z}^{\otimes \kappa - \sigma} \otimes \mathbf{I}_m \otimes  \bm{z}^{\otimes \sigma-1} \Big] \, \mathbf{f}_{p - \kappa+1} \, \bm{z}^{\otimes p - \kappa+1} \, .
\end{align*}

\sloppy This expression can be simplified further by accumulating all terms depending on $\bm{z}$ on the right. For that, it is useful to know the mixed-product property \mbox{$(\mathbf{A} \mathbf{C}) \otimes (\mathbf{B} \mathbf{D}) = (\mathbf{A} \otimes \mathbf{B}) (\mathbf{C} \otimes \mathbf{D})$} for matrices $\mathbf{A}$, $\mathbf{B}$, $\mathbf{C}$ and $\mathbf{D}$ with shapes compatible with the matrix products $\mathbf{A} \mathbf{C}$ and $\mathbf{B} \mathbf{D}$. Applying the rule to the expression inside the inner sum yields
\begin{align*}
    & \Big( \bm{z}^{\otimes \kappa - \sigma} \otimes \mathbf{I}_m \otimes  \bm{z}^{\otimes \sigma-1} \Big) \, \mathbf{f}_{p - \kappa+1} \, \bm{z}^{\otimes p - \kappa+1} = \\
    =& \Big( \mathbf{I}_{m^{\kappa - \sigma}} \, \bm{z}^{\otimes \kappa - \sigma} \Big) \otimes \Big( \mathbf{f}_{p - \kappa+1} \, \bm{z}^{\otimes p - \kappa+1} \Big) \otimes  \Big( \mathbf{I}_{m^{\sigma-1}} \, \bm{z}^{\otimes \sigma-1} \Big) = \\
    =& \Big( \mathbf{I}_{m^{\kappa - \sigma}} \otimes \mathbf{f}_{p - \kappa+1} \otimes  \mathbf{I}_{m^{\sigma-1}} \Big) \Big( \bm{z}^{\otimes \kappa - \sigma} \otimes \bm{z}^{\otimes p - \kappa+1} \otimes \bm{z}^{\otimes \sigma-1} \Big) =\\ 
    =& \Big( \mathbf{I}_{m^{\kappa - \sigma}} \otimes \mathbf{f}_{p - \kappa+1} \otimes  \mathbf{I}_{m^{\sigma-1}} \Big) \, \bm{z}^{\otimes p} \, .
\end{align*}
Above, $\mathbf{I}_a$ represents the \mbox{$a\times a$} identity matrix. Hence, the LHS at order $p$ is
\begin{align}
    & \text{LHS}\rvert_p = \sum_{\kappa = 1}^{p} \mathbf{W}_{\! \kappa} \, \mathbf{\Gamma}_{\!p, \kappa} \, \bm{z}^{\otimes p} \qquad
    & \mathbf{\Gamma}_{\!p, \kappa} := \sum_{\sigma=1}^{\kappa} \Big( \mathbf{I}_{m^{\kappa - \sigma}} \otimes \mathbf{f}_{p-\kappa+1} \otimes  \mathbf{I}_{m^{\sigma-1}} \Big)
\end{align}

\vspace{\baselineskip}\noindent
\textbf{Right-hand side}\\
The RHS expands to a sum of terms of the form
\begin{align*}
    \mathbf{A}_\kappa \, \big( \mathbf{W}_{\!\sigma_1} \bm{z}^{\otimes \sigma_1} \big) \otimes \dotsm \otimes \big( \mathbf{W}_{\!\sigma_\kappa} \bm{z}^{\otimes \sigma_\kappa} \big)
\end{align*}
Order $p$ terms must have \mbox{$\sigma_1 + \dotsm + \sigma_\kappa = p$} with $\kappa \in \mathbb{N}$ and $\bm{\sigma} \in \mathbb{N}^\kappa$. Furthermore, $\kappa$ is upper bounded, since \mbox{$\kappa = 1+\dotsm+1 \le \sigma_1 + \dotsm + \sigma_\kappa = p$}. Hence,
\begin{align*}
    \text{RHS}\rvert_p = \sum_{\kappa=1}^p \!\!\!\!\!\! \sum_{\hspace{4mm} \bm{\sigma} \in \mathbb{A}_{p, \kappa}} \hspace{-4mm} \mathbf{A}_\kappa \, \big( \mathbf{W}_{\!\sigma_1} \bm{z}^{\otimes \sigma_1} \big) \otimes \dotsm \otimes \big( \mathbf{W}_{\!\sigma_\kappa} \bm{z}^{\otimes \sigma_\kappa} \big)
\end{align*}
where \mbox{$\mathbb{A}_{p, \kappa} := \{ \bm{\sigma} \in \mathbb{N}^\kappa : \sigma_1 + \dotsm + \sigma_\kappa = p\}$}. By applying the mixed-product property, one collects all the terms depending on $\bm{z}$ on the right:
\begin{align}
    & \text{RHS}\rvert_p = \sum_{\kappa=1}^{p} \mathbf{A}_\kappa \,\mathbf{\Xi}_{p, \kappa} \bm{z}^{\otimes p} \qquad
    & \mathbf{\Xi}_{p, \kappa} := \hspace{-4mm}\sum_{\hspace{4mm} \bm{\sigma} \in \mathbb{A}_{p, \kappa}} \hspace{-4mm} \big( \mathbf{W}_{\!\sigma_1} \otimes \dotsm \otimes \mathbf{W}_{\!\sigma_\kappa} \big) \, .
\end{align}

\subsection{Order-$p$ homological equation}
Equating \mbox{$\text{LHS}\rvert_p = \text{RHS}\rvert_p$} yields
\begin{align}\label{eq:homological_base}
    \sum_{\kappa=1}^{p} \mathbf{W}_{\!\kappa} \,\mathbf{\Gamma}_{\!p, \kappa} = \sum_{\kappa=1}^{p} \mathbf{A}_\kappa \,\mathbf{\Xi}_{p, \kappa} \, .
\end{align}

Observe that $\mathbf{\Gamma}_{\!p, \kappa}$ 
is a \mbox{$m^{\kappa} \times m^{p}$} matrix whose entries are sums of $\kappa$ elements of $\mathbf{f}_{p-\kappa+1}$. Additionally, $\mathbf{\Xi}_{p, \kappa}$ 
is a \mbox{$n^\kappa \times m^p$} complex matrix whose entries are sums of products involving $\kappa$ elements of $\mathbf{W}_1$, $\mathbf{W}_2$, ... and $\mathbf{W}_{p-\kappa+1}$. Additionally, $\mathbf{\Gamma}_{\!p, \kappa}$ and $\mathbf{\Xi}_{p, \kappa}$ can be computed with the following recurrence relations
\begin{align*}
    &\mathbf{\Gamma}_{\!p, 1} = \mathbf{f}_p
    &&\mathbf{\Gamma}_{\!p+1, \kappa+1} = \mathbf{I}_{m} \otimes\mathbf{\Gamma}_{\!p, \kappa} + \mathbf{f}_{p-\kappa+1} \otimes \mathbf{I}_{m^\kappa} \, , \\
    &\mathbf{\Xi}_{p, 1} = \mathbf{W}_{\! p} 
    &&\mathbf{\Xi}_{p, \kappa+1} = \sum_{\sigma = 1}^{p-\kappa} \mathbf{W}_{\! \sigma} \otimes\mathbf{\Xi}_{p-\sigma, \kappa} \, .
\end{align*}
%

\subsection{Homological equations for each monomial}
Let \mbox{$[m] = \{1, \dotsc, m\}$}. Define an order $p$ monomial as \mbox{$\bm{z}^{\bm{\nu}} = z_{\nu_1} \dotsm z_{\nu_p}$}, where \mbox{$\bm{\nu} \in [m]^p$}. Therefore, $\bm{z}^{\otimes p}$ decomposes as
\begin{align*}
    \bm{z}^{\otimes p} = \!\!\sum_{\bm{\nu} \in [m]^p} \!\! \mathbf{e}_{\bm{\nu}} \, \bm{z}^{{\bm{\nu}}} \quad \text{with} \quad
    \mathbf{e}_{\bm{\nu}} = \mathbf{e}_{\nu_1} \otimes \mathbf{e}_{\nu_2} \otimes \dotsc \otimes \mathbf{e}_{\nu_p} \, ,
\end{align*}
where $\mathbf{e}_j$ is the $j$-th standard unit vector of dimension $m$, for \mbox{$j \in [m]$}. Similarly, $\mathbf{W}_{\!p}$ and $\mathbf{f}_p$ decompose into $m^p$ columns, each corresponding to a monomial:
\begin{align*}
    \mathbf{W}_{\!p} = \!\! \sum_{\bm{\nu} \in [m]^p} \!\! \mathbf{e}_{\bm{\nu}}^{\mathrm{T}} \otimes \mathbf{W}_{\!\bm{\nu}} \qquad
    \mathbf{f}_p = \!\! \sum_{\bm{\nu} \in [m]^p} \!\! \mathbf{e}_{\bm{\nu}}^{\mathrm{T}} \otimes \mathbf{f}_{{\bm{\nu}}} \, ,
\end{align*}
where \mbox{$\mathbf{W}_{\!\bm{\nu}} := \mathbf{W}_{\!p} \, \mathbf{e}_{\bm{\nu}}$} has shape \mbox{$n\times 1$}, and \mbox{$\mathbf{f}_{\bm{\nu}} := \mathbf{f}_{p} \, \mathbf{e}_{\bm{\nu}}$} has shape \mbox{$m\times 1$}. 

To decompose $\mathbf{\Gamma}_{\!p, \kappa}$ and $\mathbf{\Xi}_{p, \kappa}$ in monomials, we must slice $\bm{\nu}$ into $\kappa$ parts. Henceforth, given natural numbers \mbox{$b>a$}, define the slice \mbox{$\bm{\nu}[a\!:\!b]$} as the subsequence \mbox{$(\bm{\nu}_{a+1}, \dotsc, \bm{\nu}_{b})$}, which itself is a \mbox{$(b-a)$}-tuple and thereby identifies a monomial of order \mbox{$b-a$}. Therefore, for $\mathbf{\Gamma}_{\!p, \kappa}$, one has
\begin{align*}
    &\mathbf{\Gamma}_{\!\bm{\nu}, \kappa} := \mathbf{\Gamma}_{\!p, \kappa} \, \mathbf{e}_{\bm{\nu}} = \sum_{\sigma=1}^{\kappa} \Big( \mathbf{e}_{\bm{\nu}[0:\kappa-\sigma]} \otimes \mathbf{f}_{\bm{\nu}[\kappa-\sigma: p-\sigma+1]} \otimes \mathbf{e}_{\bm{\nu}[p-\sigma+1:p]}\Big) \\
    &\mathbf{\Gamma}_{\!\bm{\nu}, \kappa+1} = \mathbf{e}_{\bm{\nu}_{1}} \otimes\mathbf{\Gamma}_{\!\bm{\nu}[1:p], \kappa} + \mathbf{f}_{\bm{\nu}[0:p-\kappa]} \otimes \mathbf{e}_{\bm{\nu}[p-\kappa:p]} \, .
\end{align*}
Moreover, for $\mathbf{\Xi}_{p, \kappa}$, the same reasoning yields
\begin{align*}
     &\mathbf{\Xi}_{\bm{\nu}, \kappa} := \mathbf{\Xi}_{p, \kappa} \, \mathbf{e}_{\bm{\nu}} = \!\!\!\!\!\! \sum_{\hspace{4mm} \bm{\sigma} \in \mathbb{A}_{p, \kappa}} \hspace{-4mm} \big( \mathbf{W}_{\!\bm{\nu}[0:\sigma_1]} \otimes \dotsm \otimes \mathbf{W}_{\!\bm{\nu}[p-\sigma_k:p]} \big) \\
     &\mathbf{\Xi}_{\bm{\nu}, \kappa+1} = \sum_{\sigma = 1}^{p-\kappa} \mathbf{W}_{\bm{\nu}[0:\sigma]} \otimes \mathbf{\Xi}_{\bm{\nu}[\sigma:p], \kappa} \, .
\end{align*}

Finally, right-multiplying the order $p$ homological equation \eqref{eq:homological_base} by \mbox{$\mathbf{e}_{\bm{\nu}}$} writes the homological equation of the monomial \mbox{$\bm{z}^{\bm{\nu}} = z_{\nu_1} \dotsm z_{\nu_p}$}:
\begin{align}
    \sum_{\kappa=1}^{p} \mathbf{W}_{\!\kappa} \, \mathbf{\Gamma}_{\!\bm{\nu}, \kappa} = \sum_{\kappa=1}^{p} \mathbf{A}_\kappa \, \mathbf{\Xi}_{\bm{\nu}, \kappa} \, .
\end{align}

For the order $1$ monomial $z_j$, the homological equation becomes an eigenproblem iff the $j$-th column of $\mathbf{f}_1$ is $\mathbf{f}_{(j)} = \lambda_j \mathbf{e}_j$. Subsequently, the $j$-th column of $\mathbf{W}_{\! 1}$ becomes the eigenvector $\mathbf{W}_{\!(j)}$ associated with the eigenvalue $\lambda_j$. Expressly,
\begin{align*}
    \mathbf{W}_{\! 1} \mathbf{f}_{(j)} = \mathbf{A}_1 \mathbf{W}_{\! (j)} \implies \lambda_j \mathbf{W}_{\! (j)} = \mathbf{A}_1 \mathbf{W}_{\! (j)} \, .
\end{align*}
Higher-order homological equations simplify significantly when the linear expansion employs eigenmodes, as $\mathbf{\Gamma}_{\!p, p}$ becomes diagonal and $\mathbf{\Gamma}_{\!\bm{\nu}, p}$ reduces to
\begin{align*}
    \mathbf{\Gamma}_{\!\bm{\nu}, p} = \mathbf{e}_{\bm{\nu}} \, [\lambda_{\nu_1} \nu_1 + \lambda_{\nu_2} \nu_2 + \dotsm + \lambda_{\nu_p} \nu_p] \, .
\end{align*}

This approach introduces redundancy in monomials; for example, \mbox{$\bm{\nu}=(1,2)$} and \mbox{$\bm{\nu}=(2,1)$} both correspond to $z_1 z_2$. 
To avoid redundant computations of $\mathbf{f}_{\bm{\nu}}$ and $\mathbf{W}_{\!\bm{\nu}}$ across all permutations of $\bm{\nu}$ (since they correspond to the same monomial), sum the associated homological equations. Then choose, for example, \mbox{$\mathbf{f}_{\bm{\nu}} = \mathbf{0}$} and \mbox{$\mathbf{W}_{\!\bm{\nu}} = \mathbf{0}$} whenever $\bm{\nu}$ contains a decreasing subsequence. This method excludes terms like \mbox{$\mathbf{f}_{(2,1)}$} and \mbox{$\mathbf{W}_{\!(1,2,1)}$}, thereby simplifying higher-order $\mathbf{\Gamma}$ and $\mathbf{\Xi}$ terms. The number of homological equations to be solved also decreases. This becomes equivalent to the multi-index approach in the works~\cite{vizzaccaro2024direct,thurnher2024nonautonomous}.

As a note, when second-order ODEs are cast to first-order, the size of the system doubles and half of the variables become trivially related to the others---for instance, like \mbox{$\dot{y_i} = y_{i+n/2}$}. This greatly simplifies the problem of SSM reduction since half of the equations of the FOM are linear.
Standard structural dynamics are natively second-order in time and benefit from this simplification. Moreover, the eigensolutions of the first-order system can be obtained from the original half-sized problem, when damping is proportional.

Finally, the parametric ROM proceeds as follows:
\begin{enumerate*}
    \item isolate the forcing and parameter variables, preselecting the corresponding modes for inclusion in the ROM;
    \item extract the relevant eigenmodes from the first-order homological equation;
    \item construct and, if possible, solve higher-order homological equations to determine the associated mapping and reduced dynamics, guided by the chosen parameterization style;
    \item terminate the process upon detecting outer resonance or reaching the specified order threshold.
\end{enumerate*}

\section{Conclusion}
This work presented a framework for parametric reduction in FEM with geometry controlled via a parameter. Geometric variations are introduced without altering material properties, stress states, or inducing artificial forces. The inverse determinant in the weak-form is expanded as a power series in the parameter, with explicit expressions derived for the zeroth and first-order terms. External forcing and parameter dependence are explicitly included in an enlarged autonomous system, which can be reduced using the direct parameterization of invariant manifolds via solving a sequence of homological equations.


%


\begin{thebibliography}{6}
%

\bibitem {cabre2003parameterization1}
Cabr{\'e}, X., Fontich, E., De la Llave, R.: The parameterization method for invariant manifolds I: manifolds associated to non-resonant subspaces. 
Indiana University Mathematics Journal 52, 283--328 (2003).

\bibitem {cabre2003parameterization2}
Cabr{\'e}, X., Fontich, E., De la Llave, R.: The parameterization method for invariant manifolds II: regularity with respect to parameters. 
Indiana University Mathematics Journal 52, 329--360 (2003).

\bibitem {cabre2005parameterization3}
Cabr{\'e}, X., Fontich, E., De la Llave, R.: The parameterization method for invariant manifolds III: overview and applications. 
J. Differ. Equ. 218(2), 444--515 (2005).

\bibitem {haro2016parameterization}
Haro, A., Canadell, M., Figueras, J.-Ll., Luque, A., Mondelo, J.-M.: The parameterization method for invariant manifolds. 
Appl. Math. Sci. 195, Springer (2016).

\bibitem {guckenheimer2013nonlinear}
Guckenheimer, J., Holmes, P.: Nonlinear oscillations, dynamical systems, and bifurcations of vector fields. 
Vol. 42, Springer Science \& Business Media (2013).

\bibitem {touze2004hardening}
Touz{\'e}, C., Thomas, O., Chaigne, A.: Hardening/softening behaviour in non-linear oscillations of structural systems using non-linear normal modes. 
J. Sound Vib. 273(1–2), 77--101 (2004).

\bibitem {haller2016nonlinear}
Haller, G., Ponsioen, S.: Nonlinear normal modes and spectral submanifolds: existence, uniqueness and use in model reduction. 
Nonlinear Dyn. 86, 1493--1534 (2016).


\bibitem {vizzaccaro2024direct}
Vizzaccaro, A., Gobat, G., Frangi, A., Touz{\'e}, C.: Direct parametrisation of invariant manifolds for non-autonomous forced systems including superharmonic resonances. 
Nonlinear Dyn. 112(8), 6255--6290 (2024).

\bibitem {jain2022compute}
Jain, S., Haller, G.: How to compute invariant manifolds and their reduced dynamics in high-dimensional finite element models. 
Nonlinear Dyn. 107(2), 1417--1450 (2022).

\bibitem {opreni2023high}
Opreni, A., Vizzaccaro, A., Touz{\'e}, C., Frangi, A.: High-order direct parametrisation of invariant manifolds for model order reduction of finite element structures: application to generic forcing terms and parametrically excited systems. 
Nonlinear Dyn. 111(6), 5401--5447 (2023).

\bibitem {grolet2024high}
Grolet, A., Vizzaccaro, A., Debeurre, M., Thomas, O.: High order invariant manifold model reduction for systems with non-polynomial non-linearities: geometrically exact finite element structures and validity limit. (2024).

\bibitem {bettini2024model}
Bettini, L., Cenedese, M., Haller, G.: Model reduction to spectral submanifolds in piecewise smooth dynamical systems. 
Int. J. Non-Linear Mech. 163, 104753 (2024).

\bibitem {marconi2021higher}
Marconi, J., Tiso, P., Quadrelli, D.E., Braghin, F.: A higher-order parametric nonlinear reduced-order model for imperfect structures using Neumann expansion. 
Nonlinear Dyn. 104(4), 3039--3063 (2021).

\bibitem {morsy2025predicting}
Morsy, A.A., Tiso, P.: Predicting the variability of the dynamics of bolted joints using polynomial chaos expansion. 
Mech. Syst. Signal Process. 224, 112008 (2025).

\bibitem {thurnher2024nonautonomous}
Thurnher, T., Haller, G., Jain, S.: Nonautonomous spectral submanifolds for model reduction of nonlinear mechanical systems under parametric resonance. 
Chaos 34(7), (2024).

\bibitem {frangi2023reduced}
  Frangi, Attilio and Colombo, Alessio and Vizzaccaro, Alessandra and Touz{\'e}, Cyril: Reduced Order Modelling of Fully Coupled Electro-Mechanical Systems Through Invariant Manifolds With Applications to Microstructures.
  International Journal for Numerical Methods in Engineering,
  e7641, (2023).






\end{thebibliography}
\end{document}